\documentclass[12pt]{amsart}

\usepackage{amsthm}
\usepackage{amsmath}
\usepackage{amsfonts} %This is to have nice letters

\newcommand{\baseRing}[1]{\ensuremath{\mathbb{#1}}}
\newcommand{\Z}{\baseRing{Z}}
\newcommand{\C}{\baseRing{C}}
\newcommand{\N}{\baseRing{N}}
\newcommand{\R}{\baseRing{R}}
\newcommand{\Q}{\baseRing{Q}}

\newtheorem{theorem}{Theorem}[section]
\newtheorem{lemma}[theorem]{Lemma}

\newtheorem{definition}[theorem]{Definition}
\newtheorem{proposition}[theorem]{Proposition}
\newtheorem{example}[theorem]{Example}
\newtheorem{construction}[theorem]{Construction}
\newtheorem{conjecture}[theorem]{Conjecture}
\newtheorem{observation}[theorem]{Observation}

\def\rank{ {\rm rank}\, }
\def\vol{ {\rm vol}\, }
\def\gl{ {\rm GL} }
\def\conv{ {\rm conv} }
\def\nsupp{ {\rm nsupp} }
\def\ini{ {\rm in}\, }
\def\ind{ {\rm ind}\, }
\def\ker{ {\rm ker}\, }
\def\coker{ {\rm coker}\, }

\title{Rank jumps in Codimension 2 $A$-Hypergeometric Systems}
\author{Laura Felicia Matusevich}
\address{University of California at Berkeley}
\email{laura@math.berkeley.edu}
\date{September 11th, 2000}

\begin{document}

\begin{abstract}
The holonomic rank of the $A$-hypergeometric
system $H_A(\beta)$ is shown to depend
on the parameter vector
$\beta$ when the underlying toric ideal $I_A$ is a 
non Cohen Macaulay codimension 2 toric ideal.
The set of exceptional parameters is 
usually infinite.
\end{abstract}

\maketitle

\section{Introduction}
\label{firstsection}

$A$-hypergeometric systems are systems of linear partial 
differential equations with polynomial coefficients. In other words,
they are
left ideals in the Weyl algebra $D$,
which is the free associative algebra with
generators $x_1, \dots, x_n, \partial_1,\dots ,\partial_n$ modulo
the relations:
\[ x_i x_j = x_j x_i\; ;\; 
\partial_i \partial_j = \partial_j \partial_i \; ; \;
\partial_i x_j = x_j \partial_i + \delta_{ij} \; ,
\;\; \forall \; 1 \leq i,j \leq n \, ,\]
where $\delta_{ij}$ is the Kronecker delta.

Given a configuration of $n$ points 
$A:= \{ a_1, \ldots, a_n \} \subset \{1\} \times \Z^{d-1}$
that spans the lattice $\Z^d$ (we also think of $A=(a_{ij})$ as a $d\times n$
integer matrix of rank $d$), and a
complex vector $\beta \in \C^d$, let $H_A(\beta)$ denote the
left ideal in the Weyl algebra generated by:
\begin{equation}
\label{toricops}
\partial^u-\partial^v \, , \;\; u,v \in \N^n \; \mbox{such that} \; 
A \cdot u=A\cdot v \, ,
\end{equation}
\begin{equation}
\label{homogeneities}
\sum_{i=1}^n a_{ij} x_j \partial_j - \beta_i \, , \;\; i=1, \ldots ,d \, .
\end{equation}
The operators (\ref{toricops}) are called {\em toric operators}, and the 
operators (\ref{homogeneities}) are called {\em homogeneities}. If we set
$\theta_i = x_i \partial_i$ and $\theta$ the 
column vector whose entries are the 
$\theta_i$, then the homogeneities are simply the coordinates of the 
vector of operators $A \cdot \theta - \beta$.

The $D$-ideal $H_A(\beta)$ is called the {\em $A$-hypergeometric system 
with parameter $\beta$}. These systems, which are the object of study of 
this article,  were first introduced and studied 
by Gel'fand, Kapranov and Zelevinsky \cite{GKZ}. Solutions 
to particular instances of $H_A(\beta)$ generalize 
the classical hypergeometric functions.

\medskip

The commutative ideal of $\C[\partial_1, \dots, \partial_n]$ generated
by the toric operators will be denoted $I_A$; it is called {\em toric ideal}
or {\em lattice ideal}.
The convex hull $\conv(A)$  of the configuration $A$ is a polytope of 
dimension $d-1$. We denote its normalized volume by $\vol(A)$.
Under these hypotheses, 
$H_A(\beta)$ is a regular holonomic $D$-ideal; its {\em holonomic rank} is,
by definition, the 
common dimension of the spaces of holomorphic solutions 
of $H_A(\beta)$ around nonsingular 
points. This number is finite.

\begin{theorem}
\label{GKZ}
If $I_A$ is Cohen Macaulay, then $\rank(H_A(\beta))=\vol(A)$ for all
parameter vectors $\beta \in \C^d$.
\end{theorem}

A proof of this result, originally due to Gel'fand, Kapranov and Zelevinsky,
can be found in \cite[Section 4.3]{SST}.
The equality in the theorem can fail if $I_A$ is not Cohen Macaulay. The
following example is thoroughly analyzed in \cite{ST}.

\begin{example} 
Let $A= \left( \begin{array}{cccc} 1 & 1 & 1 & 1 \\ 0 & 1 
& 3 & 4 \end{array} \right)$. Then $\vol(A)=4$, but if we set 
$\beta=\binom{1}{2}$, we have
$\rank(H_A(\beta))=5$.
\end{example}

However, the rank of $H_A(\beta)$ is almost everywhere
equal to $\vol(A)$, as the following result shows (see \cite{A}, 
\cite[Theorem 3.5.1, Equation 4.3]{SST}).

\begin{theorem}
\label{inequality}
If $\beta$ is generic, then $\rank(H_A(\beta))=\vol(A)$.
The inequality $\rank(A) \geq \vol(A)$ always holds.
\end{theorem}

\begin{definition}
The {\bf exceptional set} of $A$ is 
\[ {\mathcal{E}}(A) = \{ \beta \in \C^d : \rank(H_A(\beta)) > \vol(A) \} \, .\]
\end{definition}
In the case $d=2$, the exceptional set is completely understood by
the following result due to Cattani, D'Andrea and Dickenstein
\cite{CDD}:

\begin{theorem}
\label{auxiliarything}
If $A = \left( \begin{array}{cccc} 1 & 1 & \ldots & 1 \\
                                   0 & a_2 & \ldots & a_n 
               \end{array} \right)$ with $0<a_2< \cdots <a_n$, then
\[ {\mathcal{E}}(A) = ( (\N A + \Z \binom{1}{0}) 
\cap (\N A + \Z \binom{1}{a_n} ) )
\backslash \N A \, . \]
This set is nonempty if and only if $I_A$ is not Cohen Macaulay.
Moreover, ${\mathcal{E}}(A)$ coincides with 
the set of parameters that maximize 
the dimension of the space of Laurent polynomial solutions of $H_A(\beta)$.
This maximum dimension is $2$.
\end{theorem}

Theorem \ref{auxiliarything} and experimental evidence give a 
basis to the following conjecture.

\begin{conjecture}
\label{conjecture}
The exceptional set ${\mathcal{E}}(A)$ of a matrix $A$ is empty if and only 
if the toric ideal $I_A$ is
Cohen Macaulay.
\end{conjecture}

\medskip

The purpose of this article is to prove Conjecture \ref{conjecture}
in the codimension $2$ case, that is, when $n-d=2$. We do this 
by explicitly constructing exceptional parameters for any
codimension $2$ non Cohen Macaulay toric ideal (see Construction
\ref{construction}). 
Our main results are:

\setcounter{section}{4}
\setcounter{theorem}{7}
\begin{theorem}
Let $\beta$ be the parameter from Construction \ref{construction}. Then 
\[ \rank(H_A(\beta)) \geq \vol(A) + 1 \, .\]
\end{theorem}

\begin{theorem}
Let $A$ be such that $I_A$ is a non Cohen Macaulay toric ideal,
with a Gale diagram whose first four rows meet each of the 
open quadrants of $\Z^2$. Let $v_1$, $v_2$, $v_4$ be as in
Construction \ref{construction}.
If $n > 4$, the $(n-4)$-dimensional affine space
parametrized by:
\[ (s_5,\dots,s_n) \longmapsto
A \cdot (v_1 e_1+v_2 e_2 -e_3 + v_4 e_4 + \sum_{i=5}^n s_i e_i)\]
is contained in the exceptional set ${\mathcal{E}}(A)$.
In particular, ${\mathcal{E}}(A)$ is an infinite set.
\end{theorem}

\setcounter{section}{1}
\setcounter{theorem}{6}

\medskip

Saito has recently announced (see \cite{Saito2}) 
that Conjecture \ref{conjecture} 
also holds when $\conv(A)$ is a simplex.

\medskip

This article is organized as follows. Section \ref{secondsection}
contains background material about 
canonical series solutions of regular holonomic systems, and in particular,
canonical $A$-hypergeometric series.
The main reference is
\cite{SST}. In
Section \ref{constructionandtools} we construct our candidates
for exceptional parameters, and develop some useful technical tools.
Section \ref{alltheproofs} contains the proofs of Theorems
\ref{rankjump} and \ref{exceptionalset}.
In Section \ref{fifthsection} we apply our methods to a concrete example,
and point out several open questions.

\section{Canonical Hypergeometric Series}
\label{secondsection}

In this section we review material concerning the series solutions of 
hypergeometric systems. We follow \cite[Sections 2.5, 3.1, 3.2, 3.4]{SST}. 
 
\begin{definition}
If $I$ is a left ideal in the Weyl algebra $D$, its {\bf distraction} 
$\tilde{I}$
is defined to be 
\[ \tilde{I}:= RI \cap \C[\theta] \, , \]
where $R= \C(x_1, \dots, x_n)\langle \partial_1,\dots,\partial_n \rangle$ 
is the
ring of differential operators with rational function coefficients, and 
$\C[\theta]=C[\theta_1,\dots,\theta_n]$ is the (commutative) subring of
$D$ generated by the operators $\theta_i=x_i \partial_i$.
\end{definition} 

The concept of distraction will allow us to define the indicial and fake
indicial ideals of a hypergeometric system. 

\begin{definition}
Let $w \in \R^n$ be a weight vector.
If $I$ is a holonomic left $D$-ideal, its {\bf indicial} ideal is 
\[ \ind_w(I) := \widetilde{\ini_{(-w,w)}(I)} = R \, \ini_{(-w,w)}(I) \cap 
\C[\theta] \, . \]
\end{definition}

The indicial ideal of a regular holonomic $D$-ideal is a zero 
dimensional ideal of the polynomial ring $\C[\theta]$. Its solutions, called
{\em exponents}, give the starting monomials 
(in a term order induced by $w$) of 
the solutions of $I$. By a {\em monomial} here we mean a product 
$x^{\alpha}=x_1^{\alpha_1}\cdots x_n^{\alpha_n}$ such that $\alpha_i \in \C$ 
and 
$x_i^{\alpha_i}=\exp(\alpha_i \log(x_i))$.

For hypergeometric ideals, there is another ideal which is closely related
to $\ind_{w}(H_A(\beta))$, but is easier to compute and understand. Its 
definition is motivated by the following facts.

\begin{proposition}{\cite[Corollary 3.1.6, Example 3.1.8]{SST}}
For \underline{generic} parameters $\beta$ we have
\[ \ind_w (H_A(\beta)) = \widetilde{\ini_w(I_A)} + \langle A \cdot \theta - 
\beta \rangle \, .\]
The containment $\supseteq$ always holds,
but $\subseteq$ can fail for non generic parameters.
\end{proposition}

The ideal $\widetilde{\ini_w(I_A)} + \langle A \cdot 
\theta - \beta \rangle$ is 
an ideal of the polynomial ring
$\C[\theta]$, called the {\em fake indicial ideal} of
$H_A(\beta)$. Its roots in affine $n$-space are called the {\em fake exponents}
of $H_A(\beta)$ with respect to $w$. Exponents are always fake exponents, and,
though the converse is not true, fake exponents are easier to describe. 
In order to do this we need to define standard pairs.
\begin{definition}
Let $M$ be a monomial ideal of $\C[\partial_1,\dots,\partial_n]$. A {\bf standard 
pair} of $M$ is a pair $(\partial^{\eta}, \sigma)$, where $\eta \in \N^n$ and
$\sigma \subset \{ 1, \dots, n \}$ subject to the following three condition:
\begin{enumerate}
\item $\eta_i = 0$ for $i \in \sigma$;
\item For all choices of integers $\mu_j \geq 0$, the monomial
$\partial^{\eta} \cdot \prod_{i \in \sigma} \partial_i^{\mu_i} $ is not in $M$.
\item For all $l \not \in \sigma$, there exist $\mu_j \geq 0$ such that
$\partial^{\eta} \cdot \partial_l^{\mu_l} \cdot 
\prod_{i \in \sigma} \partial_i^{\mu_i} $
is in $M$.
\end{enumerate}
The set of standard pairs of $M$ is denoted $S(M)$.
\end{definition}

Now we can describe the radical of the fake indicial ideal, and therefore, 
the fake exponents.

\begin{lemma}{\cite[Lemma 4.1.3]{SST}}
For any parameter vector $\beta$ and weight vector $w$ such that
$\ini_w(I_A)$ is a monomial ideal, the radical of the fake indicial ideal is
zero dimensional and equals the following intersection
\[ \bigcap_{(\partial^{\eta},\sigma) \in S(\ini_w(I_A))} 
(\langle \theta_i-\eta_i : i \not \in \sigma \rangle + \langle A \cdot \theta - \beta 
\rangle ) \, . \]
\end{lemma}

This means that, in order to compute the fake exponents, one needs only compute
the standard pairs of $\ini_w(I_A)$, and then do linear algebra.
Given a standard pair $(\partial^{\eta}, \sigma)$, the vector $v \in \C^n$
such that $v_i = \eta_i \, ,\;\; i \not \in \sigma$ and $A \cdot v = \beta$ is 
called the fake exponent with respect to $w$ corresponding to that standard pair. 
If $v$ exists, it is unique.

\medskip

Since $H_A(\beta)$ is a regular holonomic ideal, we can find a  
basis of canonical solutions of $H_A(\beta)$ with respect 
to a weight vector $w$
(see \cite[Section 2.5]{SST}). 
The elements of that basis are logarithmic series 
of the form:
\[ x^v \sum c_{v',\gamma} x^{v'} \log(x)^{\gamma}\, , \]
where $v$ is an exponent, $v' \in \ker_{\Z}(A)$, $c_{v',\gamma} \in \C$,
$\gamma \in \{0,1,\dots ,\nu-1\}^n$, and $\nu = \rank(H_A(\beta))$.

Our goal now is to describe more explicitly a basis of the space
of logarithm-free solutions of $H_A(\beta)$. The elements of this 
basis will also be canonical series.

Let $v$ be any vector in $\R^n$. Its {\em negative support} $\nsupp(v)$ is 
defined by:
\[ \nsupp(v)=\{ i \in \{1,\ldots,n\} : v_i \in \Z_{<0} \} \, . \]
The vector $v$ is said to have {\em minimum negative support} if
\begin{eqnarray*}
& v' \in \ker_{\Z}(A) \; \mbox{and} \; \nsupp(v+v') \subseteq 
  \nsupp(v) \Rightarrow &  \\
 &\Rightarrow \nsupp(v+v')=\nsupp(v) \, . &
\end{eqnarray*}
In that case, let
\[ N_v = \{ v' \in \ker_{\Z}(A) : \nsupp(v+v')=\nsupp(v) \} \, , \]
and define the following formal power series:
\begin{equation}
\label{canonicalseries}
\phi_v = \sum_{v' \in N_v} \frac{[v]_{v'_-}}{[v'+v]_{v'_-}} x^{v+v'} \;\; ,
\end{equation}
where
\[ [v]_{v'_-} = \prod_{i:v'_i<0} \prod_{j=1}^{-v'_i} (v_i-j+1) \;\;\;\; \mbox{and} 
\;\;\;\;
 [v'+v]_{v'_-} = \prod_{i:v'_i>0} \prod_{j=1}^{v'_i} (v_i+j) \, . \]

\begin{theorem}{\cite[Theorem 3.4.14, Corollary 3.4.15]{SST}}
Let $v \in \C^n$ be a fake exponent of $H_A(\beta)$ with minimum negative
support. Then $v$ is an exponent and the series $\phi_v$ defined in 
(\ref{canonicalseries}) is a canonical solution of the $A$-hypergeometric
system $H_A(\beta)$. In particular, $\phi_v$ converges in a region of $\C^n$. 
The set:
\[ \{ \phi_v : v \; \mbox{is a fake exponent with minimum negative support} \, \} \]
is a basis of the space of logarithm-free solutions of $H_A(\beta)$.
\end{theorem}

\section{Construction of exceptional parameters in codimension $2$}
\label{constructionandtools}

We assume from now on that $n-d=2$. Then $\ker_{\Z}(A)$, the integer kernel 
of $A$, is a $2$-dimensional sublattice of $\Z^n$. Let $\{ B_1, B_2 \}$
be a $\Z$-basis of $\ker_{\Z}(A)$. We think of the $B_i$ as columns
of an $n \times 2$ integer matrix $B=(b_{ji})$. The rows of $B$ form a
configuration of $n$ points in $\Z^2$. This configuration is called a
{\em Gale diagram} of $A$, and it is unique up to the action of
$\gl_2(\Z)$. The following result is contained in \cite{PS}.

\begin{theorem}
\label{fullgalecondition}
A toric ideal $I_A$ is \underline{not} Cohen Macaulay if and only if
$A$ has a Gale diagram that meets the four open quadrants of $\Z^2$.
\end{theorem}

The goal of this section is to construct exceptional parameters for $A$
when $I_A$ is a non Cohen Macaulay (codimension $2$) toric ideal.
In what follows $I_A$ will always denote such an ideal, with a Gale diagram 
$B=(b_{ij})$ that meets the four open quadrants of $\Z^2$. By 
interchanging columns of $A$ (and the corresponding rows of $B$) we may assume
that the first four rows of $B$ meet each of the four open quadrants of 
$\Z^2$, that is, $B$ is of the following form:
\[ B = \left( \begin{array}{cc} + & + \\ - & + \\ - & - \\ + & - \\ 
\vdots & \vdots \end{array} \right) \, .\]
We will need an extra assumption, that will only be used in Lemma \ref{littlebug}.
If the second and fourth row of $B$ are linearly independent, we will
assume that the cone $\{ z \in \R^2 : (B\cdot z)_2 \geq 0 \; , \; (B\cdot z)_4 \geq 0 \}$
is contained in the first quadrant. This is possible since, if this cone is not
contained in the first quadrant, it will be contained in the third.
In this case replace $B$ by $-B$. Interchanging the necessary rows of the new $B$,
we obtain a configuration as we want it.

\medskip

In the sequel, it will be very useful to compute canonical series solutions 
with respect to the weight vector $-e_3$. However, this cannot be
done if $\ini_{-e_3}(I_A)$ is not a monomial ideal. To solve this
problem while keeping all the good properties of $-e_3$ as a weight vector,
we will use the following lemma.

\begin{lemma}
\label{epsilon}
There exists $\epsilon_0 > 0$ and a generic vector $w \in \R^n$
such that, for $ 0 < \epsilon < \epsilon_0 $,
the ideal $\ini_{-e_3+\epsilon w}(I_A)=\ini_w(\ini_{-e_3}(I_A))$ 
is a monomial ideal,
and all standard pairs of $\ini_{-e_3+\epsilon  w}(I_A)$
are of the form $(\partial^{\eta}, 
\sigma = \{ 3 \}  \cup \tau)$.
\end{lemma}

\begin{proof}
The first assertion is proved using \cite[Proposition 1.13]{GBCP}
and the fact that the full dimensional cones of the Gr\"{o}bner fan of $I_A$
are exactly the cones corresponding to monomial initial ideals of $I_A$.
The second assertion is easily
proved by noticing that $\ini_{-e_3+\epsilon w}(I_A)$ is a monomial ideal
none of whose generators contain the variable $\partial_3$.
\end{proof}

Now we are ready to construct our candidates for exceptional parameters.

\begin{construction}
\label{construction}
Pick non rational
numbers $\alpha_5, \dots , \alpha_n \in \C$ such that $ \alpha_i \not
\in \Q(\alpha_5,\dots,\alpha_{i-1},\alpha_{i+1},\dots,\alpha_n)$
for $ 5 \leq i \leq n$. Remember $B_1$ and $B_2$ are the columns of $B$.
Let
\[ v = B_{1+} + B_{2+} -e_1 -e_2 - e_4 + \sum_{i=5}^n \alpha_i e_i \;\; , \]
\[ \beta = A \cdot (v-e_3) = A\cdot (v-e_3-B_1)
          = A\cdot (v-e_3-B_2) =A\cdot (v-e_3-B_1-B_2) . \]

\noindent Here we denote by $u_+$ the vector such that
$(u_+)_l= u_l$ if $u_l \geq 0$, or $(u_+)_l= 0$ otherwise, where
$u \in \Z^n$. The vector $u_-$ is defined so that $u=u_+-u_-$.
Notice that $\nsupp(v-e_3)= \{ 3 \}$, $\nsupp(v-e_3-B_1)=\{ 4 \}$,
$\nsupp(v-e_3-B_2)=\{2\}$, and $\nsupp(v-e_3-B_1-B_2)=\{1\}$.
Further, $(v-e_3)_3=(v-3_3-B_1)_4=(v-e_3-B_2)_2=(v-e_3-B_1-B_2)_1=-1$.
\end{construction}

\medskip

\begin{lemma}
\label{pseudolaurentsolutions}
The vectors $v-e_3$, $v-e_3-B_1$, $v-e_3-B_2$ and $v-e_3-B_1-B_2$
have minimum negative support. In particular, if $n=4$, $\beta \not
\in \N A$. It follows that
 $H_A(\beta)$ has
four logarithm-free solutions $\phi_{v-e_3}$, $\phi_{v-e_3-B_1}$,
$\phi_{v-e_3-B_2}$ and $\phi_{v-e_3-B_1-B_2}$, which,
for convenience in the notation, we  call
$\phi_3$, $\phi_4$, $\phi_2$ and $\phi_1$ respectively. Here
the subscripts refer to the corresponding negative supports.
\end{lemma}

\begin{proof}
First suppose that $v-e_3$ does not have minimum negative support.
Then there is $z \in \Z^2$ such that $\nsupp(v-e_3-B\cdot z)$
is strictly contained in $\nsupp(v-e_3)$. This means that
$(B\cdot z)_i \leq v_i$ for $i=1,2,4$, and $(B\cdot z)_3 < 0$.
Say $z=(z_1,z_2)^t$. If $z_1, z_2 > 0$, then $(B\cdot z)_1 \leq v_1$
does not hold. If $z_1 \leq 0, z_2 >0$, then $(B\cdot z)_2 \leq v_2$
does not hold. If $z_1 \leq 0, z_2 \leq 0 $, then $(B\cdot z)_3 < 0$
does not hold, and if $z_1>0, z_2 \leq 0$ then $(B\cdot z)_4 \leq v_4$
does not hold. All of this means that such a $z \in \Z^2$ cannot exist, 
and thus $v-e_3$ has minimum negative support. 

We show that $v-e_3-B_1$ has minimum negative support by contradiction.
Assume it does not have minimum negative support.
Then there is $z \in \Z^2$ such that $\nsupp(v-e_3-B_1-B\cdot z) = \emptyset$.
But then the negative support of $v-e_3-B \cdot (z+ (1,0)^t)$ is strictly
contained in the negative support of $v-e_3$, a contradiction. The
proofs for the other two vectors are similar. 
\end{proof}

We have found some exponents with minimum negative support of
$H_A(\beta)$. Our construction also gives an exponent with minimum 
negative support for $H_A(A\cdot v)$.

\begin{lemma}
There is only one vector with minimum negative support in the set
$\{ v + B\cdot z : z \in \Z^2 \}$. This vector is $v$, and the
corresponding logarithm-free solution of $H_A(A\cdot v)$ 
is $\phi_v=x^v$.
\end{lemma}

\begin{proof}
This follows from the same arguments that proved Lemma 
\ref{pseudolaurentsolutions}.
\end{proof}

Another interesting fact is that our construction provides an
embedded standard pair for $\ini_{-e_3+\epsilon w}(I_A)$.

\begin{lemma}
The pair
\[ ( \partial_1^{v_1}\partial_2^{v_2}\partial_4^{v_4}, \{3,5,6,\dots,n\} ) \]
is a standard pair of $\ini_{-e_3}(I_A)$, if this is a monomial
ideal. It is a standard pair of $\ini_{-e_3 + \epsilon w} (I_A)$
otherwise. Here $\epsilon$ and $w$ come from Lemma \ref{epsilon}.
\end{lemma}

\begin{proof}
To see that our candidate for standard pair satisfies
the criterion of Theorem 2.5 in \cite{hostenthomas}, we have
to show that the only integer point in a certain polytope is the
origin. This follows exactly from the same arguments of Lemma
\ref{pseudolaurentsolutions} if $\ini_{-e_3}(I_A)$ is a
monomial ideal. Otherwise, we shrink $\epsilon$ so that the same
arguments will work when we use the weight vector $-e_3+\epsilon w$.

We also need to find elements in $\Z^2$ that belong to that
polytope when one of the defining inequalities is removed. 
Those elements will be $(1,0)^t$, $(0,1)^t$ and $(1,1)^t$.
\end{proof}

We want to show that $\rank(H_A(\beta)) > \vol(A)$. In view of
Theorem \ref{inequality}, one way to do this is to show
that $\rank(H_A(\beta))$ is strictly greater than $\rank(H_A(A\cdot v))$.
In order to compare this two numbers, we need a link between $H_A(\beta)$
and $H_A(A\cdot v)$. This is provided by the following
$D$-module map (see \cite[Section 4.5]{SST}):
\[ \partial_3 : D/H_A(\beta) \longrightarrow D/H_A(A\cdot v) \; .\]
This $D$-module map induces a vector space homomorphism in the opposite 
direction between the solution spaces of the corresponding hypergeometric ideals,
namely, if $\varphi$ is a solution of $H_A(A\cdot v)$, then $\partial_3 \varphi$
is a solution of $H_A(\beta)$. 

Our strategy to show that $\rank(H_A(\beta)) > \rank (H_A(A\cdot v))$ will
be as follows. First, characterize the kernel of $\partial_3$ (as a map between
solution spaces). There is an obvious element of this kernel, namely
the function $\phi_v=x^v$. After we have done that, we will construct, for
each element of a vector space basis of $\ker(\partial_3)$, a nonzero
function in the cokernel of $\partial_3$. However, for the function $\phi_v$ 
(which will belong to that generating set) we will construct at least two functions
in $\coker(\partial_3)$. After showing all of the functions thus constructed are 
linearly independent, we will conclude $\dim(\coker(\partial_3))
\geq \dim(\ker(\partial_3))+1$. This will imply the desired result 
(that is, that $\rank(H_A(\beta)) > \rank (H_A(A\cdot v)$)
using elementary linear algebra.

Before we can look at the kernel and cokernel of $\partial_3$, we need
a couple of technical facts.

\medskip

\begin{observation}
\label{keyobservation}
{\rm 
Let $\psi$ be a solution of $H_A(A\cdot v)$.
This function is of the form:
\[ \psi = \sum 
c_{\alpha,\gamma}x^{\alpha} \log(x)^{\gamma} \, ,\]
where $A\alpha=A\cdot v$, 
$\nu = \rank(H_A(A\cdot v))$, $\gamma \in \{0,1,\dots ,\nu-1\}^n$,
and $\log(x)^{\gamma}=\log(x_1)^{\gamma_1}\cdots \log(x_n)^{\gamma_n}$.

The set ${\mathcal{S}}:=\{ \gamma \in [0,\nu-1]^n\cap\N^n : \exists 
\alpha \in \C^n \,
\mbox{such that} \; c_{\alpha,\gamma} \neq 0 \}$ 
is partially ordered with respect to:
\[ (\gamma_1,\ldots,\gamma_n) \leq (\gamma_1',\ldots,\gamma_n') 
\Leftrightarrow \gamma_i \leq \gamma_i', \;
i=1,\ldots , n \, .\]

Denote by ${\mathcal{S}}_{\max}$ the set of maximal elements of 
${\mathcal{S}}$.
Let $\delta \in {\mathcal{S}}_{\max}$ and  $f = \sum_{\alpha \in \C^n}
c_{\alpha,\delta}x^{\alpha}$. Write
\[ \psi = \psi_{\delta} + \log(x)^{\delta} \, f \, ,\]
so that the logarithmic terms in $\psi_{\delta}$ are either less than
or incomparable to $\delta$.
If $P$ is a differential operator that annihilates $\psi$, we have:
\[0 = P \psi = P\psi_{\delta}+\log(x)^{\delta} Pf+
\mbox{terms whose $\log$ factor is lower than
$\delta$} .\]
Since $P\psi_{\delta}$ is a sum of terms whose $\log$ factor is 
either lower than $\delta$
or incomparable to $\delta$, we conclude that $Pf$ must be zero. 
This implies that
$f$ is a logarithm-free $A$-hypergeometric function of degree $A\cdot v$. }
Moreover, if $\partial_3 \psi$ is logarithm-free, then
$\partial_3 f$ must vanish.
\end{observation}

The following lemma is used to analyze
the kernel and cokernel of the map $\partial_3$, and it will be used 
repeatedly in the sequel. Its proof is inspired by the proofs of
Theorems 2.5 and 3.1 in \cite{hostenthomas}.

\begin{lemma}
\label{bigbug}
Let $u$ be a fake exponent of $H_A(A \cdot v)$ such that $u_3=0$,
corresponding to a 
standard pair
$(\partial^{\eta},\sigma=\{3\} \cup \tau)$ of 
$\ini_{-e_3+\epsilon w}(I_A)$. 
Here $\epsilon$ and $w$ are chosen so that $\ini_{-e_3+\epsilon' w}
=\ini_w(\ini_{-e_3}(I_A))$ for all $0<\epsilon' \leq \epsilon$ and this 
is a monomial ideal.
Then there exists a set 
${\mathcal{I}} \subseteq \{1,2,4,\dots,n \} \backslash \tau$ of cardinality
$2$ such that,  for each
$i \in {\mathcal{I}}$, we can find a vector $z^{(i)}
\in \Z^2$ that satisfies the following three properties:
\begin{enumerate}
\item $(B\cdot z^{(i)})_i > \eta_i$,
\item $(B\cdot z^{(i)})_j \leq u_j \;\;\;$ for all $j \neq 3,i$ such that $u_j \in \N\, ,$
\item $(B\cdot z^{(i)})_3 < 0$.
\end{enumerate}
Moreover, ${\mathcal{I}}$ can be chosen so that the rows of $B$
indexed by ${\mathcal{I}}$ are linearly independent.
\end{lemma}

\begin{proof}

Fix $l \not \in \sigma = \{3\} \cup \tau$.
Let  $\mu \in \N^n$ such that $\mu_j = u_j$ if $u_j \in \N$; $\mu_j = 0$
otherwise. Observe that $\mu_j = \eta_j$ for $j \not \in \sigma$.
For $v' \in \N^n$ we define, following \cite{hostenthomas}:
\[ P_{v'}(0) = \{ y \in \R^2 : B\cdot y \leq v'; 
-(-e_3+\epsilon w)^t (B\cdot y) \leq 0 \} \, .\]
Following the proof of Theorem 2.5 in \cite{hostenthomas}, we see that,
for $l \not \in \sigma$ there exists a positive integer $m_l$ such
that $P_{\mu+m_l e_l}(0)$ contains a nonzero integer
vector $z^{(l)}\in \Z^2$. 
It must satisfy $-(-e_3+\epsilon w)^t (B\cdot z^{(l)}) <0$. 
The reason for this is that, since
$\ini_{-e_3+ \epsilon w}(I_A)$ is a monomial ideal, 
there exists a unique solution
of the integer program
\[ \mbox{ minimize} \,    -(-e_3+\epsilon w)^t (B\cdot z) \; 
\mbox{subject to} \;
z \in P_{\mu + m_l e_l}
\cap \Z^2 \; , \]
where $P_{\mu + m_l e_l}:=\{ y \in \R^2 : B\cdot y \leq \mu + m_l e_l \}$ (see 
\cite[Section 2]{hostenthomas}), and we can choose $z^{(l)}$ as that solution.

The vectors $z^{(l)}$ are almost what we want, except that we cannot
a priori guarantee that $(B\cdot z^{(l)})_3 < 0$, even if we look at
all the polytopes $P_{\mu+m e_l}(0)$ for $m \geq m_l$, that is,
even if we look at the (possibly unbounded) polyhedron:
\[ R^l:= \{ y \in \R^2 : (B \cdot y)_j \leq \mu_j,  \; j \neq 
l; -(-e_3+\epsilon w)^t (B\cdot y) \leq 0 \} \, .\]
However, we may assume $(B\cdot z^{(l)})_3 \leq 0$ since
we can always choose $\epsilon$ small enough so that 
a feasible point that satisfies $(B\cdot z)_3\leq 0$ is better
than one that satisfies $(B\cdot z)_3 > 0$.

The following notation is very convenient:
\[ P^{\bar{\sigma}}_{\eta}(0) := \{
z \in \R^2 : (B\cdot z)_i \leq \eta_i \; , i \not \in \sigma 
; -(-e_3+\epsilon w)^t (B\cdot z) \leq 0\} \, ,\]
\[ E^l:=\left\{ \begin{array}{rl} z \in \R^2 :& \begin{array}{l}
(B \cdot z)_j \leq \mu_j,  \;j \neq 
l; (B\cdot z)_l> \eta_l; \\  -(-e_3+\epsilon w)^t (B\cdot z) \leq 0 
\end{array} \end{array}
\right\}\, . \]
Notice that $E^l=R^l \backslash P^{\bar{\sigma}}_{\eta}(0)$.

Let us first deal with the case when the hyperplane $\{ (B\cdot z)_l =0 \}$
is parallel to $\{ (B \cdot z)_3 = 0\}$, that is, when there exists
$\lambda \in \Q$ such that $e_l^t B = \lambda e_3^t B$. 
We know that $u = v - B\cdot y$ for some $y \in \C^2$. Since $u_3=0=v_3$,
we have $(B\cdot y)_3=0$ which implies $(B\cdot y)_l=0$, so that
$u_l=v_l$.
But $u_l$ is an integer. By construction of $v$, this implies
$l=1$ and $\lambda <0$ (remember that the only integer coordinates
of $v$ are the first four).
Now $v_1 < (B\cdot z^{(1)})_1= \lambda (B\cdot z^{(1)})_3$ and
$\lambda<0$ imply that $(B\cdot z^{(1)})_3< v_l/\lambda <0$.

Now fix $l \not \in \sigma$ such that the $l$-th row of $B$ is not a multiple
of the third one, and suppose that the integer program 
\[ \mbox{minimize} \, -(-e_3+\epsilon w)^t(B\cdot z) 
\, \mbox{subject to} \;  z \in R^l\cap \Z^2  \]
is unbounded, and every bounded subprogram has its solution
on the hyperplane $\{ (B\cdot z)_3 =0 \}$. Then 
$R^l\cap \Z^2 \cap \{ (B\cdot z)_3 = 0 \}$ is an infinite set. 
Notice that $R^l$ is not contained in the half-space $\{ (B\cdot z)_3 
\geq 0\}$,
since the defining inequalities of $R^l$ given by rows that are multiples
of the first row of $B$ are of the form $(B\cdot z)_3 \leq 0$.
This follows from similar arguments as those in the
preceeding paragraph. But now the set 
$R^l \cap \{ (B\cdot z)_3 \leq 0\}$ contains infinitely many lattice
points on the hyperplane $\{ (B\cdot z)_3=0\}$, is not itself 
contained in this hyperplane, but is a subset of 
$\{ z \in \R^2 : -1 < (B\cdot z)_3 \leq 0 \}$. This is impossible.

Thus, if $z^{(l)}$ satisfies $(B\cdot z)_3=0$, the integer program:
\[ \mbox{minimize} \, -(-e_3+\epsilon w)^t(B\cdot z) \,
\mbox{subject to} \; z \in R^l\cap \Z^2  \]
must be bounded. 
Let ${\mathcal{J}}\subseteq \{1,2,4,\dots ,n\} \backslash
\tau$ be the set of all such indices $l$, with $z^{(l)}$ the 
(unique) solution to the corresponding integer program.
We can now follow the proof of Theorem 3.1 in \cite{hostenthomas}
to show that $\langle \partial_i : i \not \in \sigma \cup {\mathcal{J}}
\rangle$ is an associated prime of $\ini_{-e_3+\epsilon w}(I_A)$.

Now let ${\mathcal{I}}$ be such that $\langle \partial_i : i \in 
{\mathcal{I}} \rangle$ is a minimal
prime of $\ini_{e_3+\epsilon w}(I_A)$ containing 
$\langle \partial_i : i \not \in \sigma \cup {\mathcal{J}} \rangle$.
Then the vectors $z^{(l)}$ for $l \in {\mathcal{I}}$ satisfy
all the desired properties, and the cardinality of ${\mathcal{I}}$
is $2$. 

The only thing we still have to show is that the rows of
$B$ indexed by ${\mathcal{I}}$ are linearly independent.
To see this, let $(\partial^{\eta'},\sigma':= \{ 1, \dots, n\} \backslash
{\mathcal{I}})$ be a standard pair of $\ini_{-e_3+\epsilon w}(I_A)$, and 
look at the set:
\[ P^{\bar{\sigma'}}_{\eta'}(0) = 
\{ z \in \R^2 : (B \cdot z)_i \leq \eta'_i\; , i \not \in \sigma' ;
-(-e_3+\epsilon w)^t (B\cdot z) \leq 0 \} \, ,\]
which, by Theorem 2.5 in \cite{hostenthomas}, is a polytope.
If the rows of $B$ indexed by ${\mathcal{I}}$ are not linearly independent,
then $2\times 2$ matrix whose rows are those rows of $B$ has a nontrivial 
kernel. Hence there exists $y \in \R^2$ such that
$(B\cdot y)_i = 0$ for all $i \not \in \sigma'$. Since
all the $\eta'_i$ are nonnegative, this means that 
$P^{\bar{\sigma'}}_{\eta'}(0)$ contains at least half of the line 
$\{ s y : s \in \R\}$, contradicting the fact
that $P^{\bar{\sigma'}}_{\eta'}(0)$ is a bounded set. This concludes the proof.
\end{proof}

Notice that a stronger result holds for the fake exponent
$v$ corresponding to the standard pair $(\partial_1^{v_1}\partial_2^{v_2}
\partial_3^{v_3},
\{3, 5,\dots ,n\})$, namely the three vectors $(1,1)^t$,$(0,1)^t$ and $(1,0)^t$
satisfy the properties required of the vectors
$z^{(l)}$ in Lemma \ref{bigbug}.

\section{The structure of the map $\partial_3$}
\label{alltheproofs}

In this section we study the kernel and cokernel of the map $\partial_3$
between the solution spaces of $H_A(A\cdot v)$ and $H_A(\beta)$. 
The following proposition is the first step towards describing its kernel.

\begin{proposition}
\label{monokernel}
If $\varphi$ is a canonical logarithm-free series 
solution of $H_A(A\cdot v)$ and 
$\partial_3\varphi$ belongs to 
$\mbox{Span}\, \{ \phi_1,\phi_2, \phi_3, \phi_4\}$, where the 
functions $\phi_i$
are the logarithm-free canonical series from
Lemma \ref{pseudolaurentsolutions},
then $\varphi$ is a monomial and $\partial_3 \varphi=0$.
\end{proposition}

\begin{proof}

We compute canonical series with respect to the weight vector $-e_3$,
as in \cite[Sections 2.5, 3.4]{SST}, assuming that $\ini_{-e_3}(I_A)$
is a monomial ideal. We will deal with the case when $\ini_{-e_3}(I_A)$
is not monomial later.

The logarithm-free canonical solutions of $H_A(A\cdot v)$ are of the form
\[ \phi_u = \sum_{v' \in N_u} \frac{[u]_{v'_-}}{[v'+u]_{v'_-}} x^{u+v'}\, , \]
where
$u$ is a fake exponent of minimum negative support. The fact that
$u$ is a fake exponent means that
there exists a standard pair $(\partial^{\eta}, \sigma)$ of $\ini_{-e_3}(I_A)$,
with $\sigma=\{3 \} \cup \tau$, such that $u$ is the unique vector
satisfying $A\cdot u=A\cdot v$ and $u_i=\eta_i, \; i \not \in \sigma$. 

The only fake exponent with minimum negative support in 
$\{ v+B\cdot z: z \in \Z^2\}$
is $v$, whose canonical solution is $x^v$, and this function satisfies 
$\partial_3 x^v = 0$.
Let $u$ be a fake exponent
with minimum negative support that does not differ with $v$
by an integer vector.
Call $\varphi$ its canonical solution. If $\partial_3 \varphi$
belongs to $\mbox{Span}\, \{ \phi_1,\phi_2, \phi_3, \phi_4\}$,
it is clear that we must have $\partial_3 \varphi=0$,
that is, $\varphi$ must be a constant function with respect to $x_3$.
In particular, we need $u_3=0$.

If $v'=B \cdot z$ is an element of $N_u$, then it must satisfy the inequalities
\[ (B\cdot z)_i \geq -\eta_i \; , \; i \not \in \sigma 
\; ; \; \; (B \cdot z)_3 \geq 0 \, .\]
But the set
\[ P_{\eta}^{\bar{\sigma}}(0) := \{ z \in \R^2 : (B\cdot z)_i \leq \eta_i 
\; , \; i \not \in \sigma \; ; \; \; (B\cdot z)_3 \leq 0 \} \]
intersects the lattice $\Z^2$ only at $0$ (see 
\cite[Theorem 2.5]{hostenthomas}). Switching 
the inequality signs, we conclude $N_u=\{ 0 \}$, so that
$\varphi=x^u$. 

Now, if $\ini_{-e_3}(I_A)$ is not a monomial ideal, take
$w$ and $\epsilon_0$ as in Lemma \ref{epsilon}. We can choose $0 < 
\epsilon < \epsilon_0$
so that the polytopes
\[ P_{\eta}^{\bar{\sigma}}(0) := \{ z \in \R^2 : (B\cdot z)_i \leq \eta_i 
\; , \; i \not \in \sigma \; ; \; \; -(-e_3+\epsilon w)^t (B\cdot z) 
\leq 0 \} \]
and
\[ \{ z \in \R^2 : (B\cdot z)_i \leq \eta_i 
\; , \; i \not \in \sigma \; ; \; \; (B\cdot z)_3 \leq 0 \} \]
have the same integer points. Now the previous reasoning applies when we
compute canonical series with respect to $-e_3+\epsilon w$ instead of $-e_3$.
\end{proof}

We are now ready to characterize the kernel of $\partial_3$ as a map
between solution spaces.

\begin{theorem}
\label{kernel}
The kernel of the map 
\[ \partial_3 : \{ \, \mbox{Solutions of} \;  H_A(A\cdot v) \} \longrightarrow
\{ \, \mbox{Solutions of} \;  H_A(\beta) \} \]
is spanned by
\[ \left\{ x^u : \begin{array}{c}
               \, \mbox{$u$ is (a fake) exponent with minimum} \\
               \mbox{negative support such that $u_3=0$} 
\end{array} \right\} \; .\]
\end{theorem}

\begin{proof}
It is clear that the functions described above belong to the kernel of 
$\partial_3$.
Suppose first that $\varphi$ is a logarithm-free solution of $H_A(A\cdot v)$
that is constant with respect to $x_3$.
We compute canonical series with respect to the weight vector $-e_3$.
If this cannot be done (that is, if $\ini_{-e_3}(I_A)$ is not a monomial
ideal) we replace this weight by $-e_3+\epsilon w$ from Lemma \ref{epsilon}
with $\epsilon$ small enough so that the ideas still work.

Now $\varphi$
is a linear combination of
logarithm-free canonical series
(with respect to the weight $-e_3$), each corresponding to a fake 
exponent with 
minimum negative support.
Say $\varphi = \sum c_{u^{(i)}} \phi_{u^{(i)}}$, where $c_{u^{(i)}} \in 
\C$ and $u^{(i)}$  are the
exponents with minimum negative support.

By taking initials, we see that at least one of those exponents must have
its first coordinate equal to zero. Call that exponent $u$. But then,
by the proof of Proposition \ref{monokernel}, 
the canonical series corresponding to $u$ is
$x^u$, and this function belongs to our candidate spanning set. 
Subtracting $c_ux^u$ to $\varphi$ and repeating the process, we conclude that 
$\varphi$ is a linear combination of the functions in our
candidate spanning set.

Our task now is to show that no logarithmic solution of $H_A(A\cdot v)$ can be 
constant with respect to $x_3$. 

Let $\psi$ be a (possibly 
logarithmic) solution of $H_A(A\cdot v)$ and suppose that 
$\partial_3 \psi =0$. The function $\psi$ is a linear combination
of canonical series.
We write $\psi=\varphi_1+\cdots+\varphi_k$
where in each $\varphi_i$ we collect all canonical series appearing as summands
in $\psi$ whose corresponding exponents differ by integer vectors.
Then there exist $u^{(i)}$ exponents with minimum negative support and
first coordinate equal to zero, such that 
$\varphi_i = \sum c_{\gamma, \alpha} x^{\alpha} \log(x)^{\gamma}$,
where $c_{\alpha,\gamma} \neq 0 \Rightarrow \alpha-u^{(i)} \in \Z^n$.
Also notice that each $\varphi_i$ must be constant with respect to $x_3$.

We must show that each function $\varphi_i$ must be logarithm-free.
Pick one of those functions $\varphi_i$ and the exponent $u^{(i)}$.
We will now drop the index $i$ for convenience in the notation.
Write $\varphi$ in the form of Observation \ref{keyobservation}. In this case
$f=x^u$ for any $\delta \in {\mathcal{S}}_{\max}$ by construction of $\varphi$.
Now we apply Lemma \ref{bigbug} to the exponent $u$. 
Let $j \in {\mathcal{I}}$, write $z$ for the vector $z^{(j)}$ and
let $\delta \in {\mathcal{S}}_{\max}$ be maximal with respect to the
$j$-th coordinate. Remember $\varphi = \varphi_{\delta}
+c_{\delta} x^u \log(x)^{\delta}$,
where $\varphi_{\delta}$ contains only terms in $\log$ that are either
less than $\delta$ or incomparable to $\delta$.
We know that $\partial^{(B\cdot z)_-} \varphi=0$, since
$(B\cdot z)_3 <0$. Then
\begin{eqnarray*}
0 & = & \partial^{(B\cdot z)_-} \varphi \\
  & = & \partial^{(B\cdot z)_+} \varphi \\
  & = & \partial^{(B\cdot z)_+} \varphi_{\delta} + \partial^{(B\cdot z)_+}
x^u \log(x)^{\delta}
\end{eqnarray*}
All the terms that come from $\partial^{(B\cdot z)_+} x^u \log(x)^{\delta}$
by applying the product rule are either zero or must be cancelled 
by something from $\partial^{(B\cdot z)_+} \varphi_{\delta}$.
As a matter of fact, 
$\partial^{(B\cdot z)_+} x^u \log(x)^{\delta}$
has a nonzero term which is a multiple of
\[  \frac{(\partial^{(B\cdot z)_+ 
+ (-(B\cdot z)_j + \eta_j)e_j} x^u)
\log(x)^{\delta-((B\cdot z)_j-\eta_j) e_j}}{x_j^{(B\cdot z)_j-\eta_j}} \]
if $(B\cdot z)_j - \eta_j \leq \delta_j$, or of
\[  \frac{(\partial^{(B\cdot z)_+ 
+ (-(B\cdot z)_j + \eta_j)e_j} x^u)
\log(x)^{\delta-\delta_j e_j}}{x_j^{(B\cdot z)_j-\eta_j}}\]
otherwise. The numerators of these fractions are nonzero by 
construction of $z$.
Then we have a 
sub-series $g$ of $\varphi_{\delta}$ such that
\[\partial^{(B\cdot z)_+} (g-x^u\log(x)^{\delta})=0 \, . \]
This means that $g-x^u\log(x)^{\delta}$ is a polynomial in the variable $x_j$,
which contradicts the fact that $\varphi_{\delta}$ contains no
term in $\log(x)^{\delta}$.
This implies that $\delta_j=0$, so that $\varphi$ contains no
$\log(x_j)$, and this is true for all $j \in {\mathcal{I}}$.

Now pick any $l \not \in {\mathcal{I}}$, and $\delta \in {\mathcal{S}}_{\max}$
maximal with respect to the $l$-th coordinate. As before,
$\varphi = \varphi_{\delta}+ c_{\delta}x^u \log(x)^{\delta}$.
Of course, since $x^u$ is itself a hypergeometric function constant
with respect to $x_3$, we may assume that $\varphi$ has no term in $x^u$.
This and the homogeneity equations (\ref{homogeneities}) imply that
there is a sub-sum of $\varphi$ of the form
$x^u \sum_{k=1}^n c_k \log(x)^{\delta-e_l+e_k}$, where
$(c_1,\dots,c_n)^t$ belongs to the kernel of $A$, and there are no other 
terms in $x^u \log(x)^{\delta-e_l+e_k}$ in $\varphi$.

From our previous reasoning, we know that $c_j=0$ for all 
$j \in {\mathcal{I}}$. Since the rows of $B$ indexed by ${\mathcal{I}}$
are linearly independent (and the columns of $B$ span the kernel of $A$),
we conclude that $(c_1,\dots,c_n)^t=0$. In particular, $c_l=c_{\delta}=0$.
This completes the proof that $\varphi$ is logarithm-free.
\end{proof}

\noindent {\em Remark:}
Currently, all the examples where we have computed the map
$\partial_3$ have a $1$-dimensional kernel. However, all these
examples are small, so we believe that
there will be examples where $\partial_3$ has a  higher-dimensional
kernel.

\medskip

We want to compute the dimension of the solution space of $H_A(\beta)$
using information about the dimension of the kernel and cokernel of the
map $\partial_3$. In particular, our goal is to show that
the sum of the dimension of the image of $\partial_3$ and the
dimension of the cokernel of $\partial_3$ is at least the dimension
of the solution space of $H_A(A\cdot v)$ plus one.
The next step in this direction is to find linearly independent
solutions of $H_A(\beta)$ not lying in the image of $\partial_3$ corresponding
to the elements of the kernel of $\partial_3$.

\begin{lemma}
\label{littlebug}
Let $u$ be a fake exponent of $H_A(A\cdot v)$ with minimum negative support 
corresponding to a standard pair $(\partial^{\eta}, 
\sigma = \{3 \} \cup \tau)$,
and assume that $u_3=0$. Then $u-e_3$ is the
fake exponent of $H_A(\beta)$ corresponding to  
$(\partial^{\eta}, \sigma = \{3 \} \cup \tau)$, 
and it has minimum negative support.
\end{lemma}

\begin{proof}
That $u-e_3$ is the fake exponent of $H_A(\beta)$ corresponding to
$(\partial^{\eta}, \sigma = \{3 \} \cup \tau)$ follows from the fact that 
$3 \in \sigma$ (and 
that we have only modified the third coordinate of $u$).

Now we have to show that $u-e_3$ has minimum negative support. 
We know that $u_i \in \N$ for $i \not \in \tau$, so that $u$ has at least 
three integer coordinates. If it has exactly those integer coordinates,
or if its integer coordinates are all greater than or equal to zero, 
then $\nsupp(u-e_3)=\{ 3 \}$. It follows that it has minimum negative
support. To see this, suppose $\nsupp(u-e_3-(B\cdot z))$ is
strictly contained in $\nsupp(u-e_3)$ for some $z \in \Z^2$. This means that
$(B\cdot z)_i \leq \eta_i$, for $i \not \in \sigma$, and
$(B\cdot z)_3 < 0$. Then $z \in P^{\bar{\sigma}}_{\eta}(0) \cap \Z^2 = \{
0 \}$, a contradiction.

Now assume that $u$ has some negative integer coordinates,
and write $u = v-B\cdot y$ for some $y \in \C^2$. Then $u$
has at least four integer coordinates. We claim that in that case,
$u$ has exactly four integer coordinates, and they are the first four.
To show this claim that we will use the numbers
$\alpha_i$ from Construction \ref{construction}. We know $u_3=0$, so
that $(B\cdot y)_3=0$. Suppose that $u_j \in \Z$ for some $j>4$.
Then the $j$-th column of $B$ and the third column of $B$ are linearly 
independent, because otherwise, we would have $(B\cdot y)_j =0$
so that $u_j=v_j \not \in \Z$. This means that $y \in \Q(\alpha_j) 
\backslash \Q$. But now the construction of the numbers $\alpha_i$ 
implies that the only integer coordinates of $u$ must be the 
third one, the $j$-th one, and maybe the first one (if the first row of $B$ is
a multiple of the third). We obtain
a contradiction. Thus, the only integer 
coordinates of $u$ are the first four. Moreover, $u$ has some negative integer 
coordinates. This can only happen if 
$(\partial^{\eta}, \sigma = \{3 \} \cup \tau)$ is a top dimensional standard
pair, $\{1,\dots,n\} \backslash \sigma$ is strictly contained
in $\{1,2,4\}$, and $u_j$ is a negative integer, where
$j$ is the only element of $\{1,2,4\} \cap \sigma$.

Assume that $u-e_3$ does not have minimum negative support, and
pick $z \in \Z^2$ such that $\nsupp(u-e_3-(B\cdot z))$ is
strictly contained in $\nsupp(u-e_3)$. Looking at 
$P^{\bar{\sigma}}_{\eta}(0)$, we conclude that
we cannot have $(B\cdot z)_3 \leq 0$. Then $(B\cdot z)_3 > 0$ and
$(B\cdot z)_j \leq u_j < 0$. It follows that $u-B\cdot z$ has minimum 
negative support $\{ 3 \}$ (and is thus an exponent of $H_A(A\cdot v)$).
We will show that $u-B\cdot z$ actually does not have
minimum negative support. This contradiction will imply the desired conclusion
about $u-e_3$.

In order to show that 
$u-B\cdot z$ does not have minimum negative support, we need to find a 
vector $\tilde{z} \in \Z^2$ such that the negative support of
$u-B\cdot z - B \cdot \tilde{z}$ is empty.
We know that $u-B\cdot z = v - B \cdot (y+z)$, $(B\cdot (y+z))_3 > 0$
and $(y+z) \neq 0$.
We have the following cases:
\begin{enumerate}
\item $(B\cdot (y+z))_1 < 0, (B\cdot (y+z))_2 \geq 0, (B\cdot (y+z))_4 < 0$ ,
\item $(B\cdot (y+z))_1 < 0, (B\cdot (y+z))_2 < 0, (B\cdot (y+z))_4 \geq 0$ ,
\item $(B\cdot (y+z))_1 < 0, (B\cdot (y+z))_2 < 0, (B\cdot (y+z))_4 < 0$ ,
\item $(B\cdot (y+z))_1 < 0, (B\cdot (y+z))_2 \geq 0, (B\cdot (y+z))_4 \geq 0$ ,
\item $(B\cdot (y+z))_1 \geq 0, (B\cdot (y+z))_2 
\geq 0, (B\cdot (y+z))_4 < 0$ ,
\item $(B\cdot (y+z))_1 \geq 0, (B\cdot (y+z))_2 
< 0, (B\cdot (y+z))_4 \geq 0$ .
\end{enumerate}
In case $1$, $\nsupp(v-B\cdot(y+z)-B_1)$ is contained in
$\nsupp(v-B\cdot(y+z))$. In case $2$,  $\nsupp(v-B\cdot(y+z)-B_2)$ 
is contained in $\nsupp(v-B\cdot(y+z))$. In case $3$,
$\nsupp(v-B\cdot(y+z)-(B_1+B_2))$ is contained in 
$\nsupp(v-B\cdot(y+z))$. In case $4$, 
$\nsupp(v-B\cdot(y+z)-(B_1+B_2))$ is contained in 
$\nsupp(v-B\cdot(y+z))$.
In case $5$
$\nsupp(v-B\cdot(y+z)-B_1)$ is contained in
$\nsupp(v-B\cdot(y+z))$. In case $6$,
$\nsupp(v-B\cdot(y+z)-B_2)$ is contained in
$\nsupp(v-B\cdot(y+z))$.
Cases $1$, $2$ and $3$ follow directly from the construction of $v$.
For case $4$, remember that we assumed at the beginning
of Section  \ref{constructionandtools} that either
the second and fourth rows of $B$ are linearly dependent, or
the cone $\{ z\in \R^2 : (B\cdot z)_2 \geq 0 \, , \, (B\cdot z)_4 \geq 0 \}$
is contained in the first quadrant.
This means that the only way case $4$ could happen is if 
the second and fourth rows of $B$ are linearly dependent, and
$(B\cdot z)_2=(B\cdot z)_4=0$. Then our assertion about
negative supports follows by direct verification.
Finally, let us do case $5$. Case $6$ will be similar.

Since $b_{21} < 0$ and $(v-B\cdot(y+z))_2\geq 0$, we have
$(v-B\cdot(y+z)-B_1)_2 = (v-B\cdot(y+z))_2 - b_{21}\geq 0$. Since
$(B\cdot (y+z))_4$ is a negative integer, we have
$(v-B\cdot(y+z)-B_1)_4 = b_{41}-1-(B\cdot(y+z))_4-b_{41} \geq 0$.
The inequalities that $y+z$ satisfies imply that this vector belongs to the
second quadrant of $\Z^2$. Its second coordinate is
strictly less than one. To see this, remember that 
$(B\cdot(y+z))_2 \leq v_2=b_{22}-1$. The line
$\{ (s_1,s_2) \in \R^2: (B\cdot(s_1,s_2)^t)_2 =v_2\}$ 
cuts the vertical axis of $\R^2$ above
zero and strictly below one (this is because the line
$\{ (s_1,s_2):(B\cdot (s_1,s_2)^t)_2 = v_2+1\}$ 
cuts the vertical axis at height $1$).
It follows that $0 \leq y_2+z_2 <1$.
We have:
\begin{eqnarray*}
(v-B\cdot(y+z)-B_1)_1 &=& b_{11}+b_{12}-1 -b_{11}(z_1+y_1)-\\
& &  - b_{12}(z_2+y_2)-b_{11} \\
 &= & -b_{11}(z_1+y_1)+b_{12}-b_{12}(z_2+y_2) -1 
\end{eqnarray*}
We know $ -b_{11}(z_1+y_1) \geq 0$, since $z_1+y_1 \leq 0$.
We also know $b_{12}-b_{12}(z_2+y_2) \geq 0$. The sum of this two
numbers is an integer (since $(v-B\cdot(y+z)-B_1)_1$ is an integer),
so it must be greater than or equal to $1$. This implies 
$(v-B\cdot(y+z)-B_1)_1 \geq 0$, and concludes the proof that 
$\nsupp(v-B\cdot(y+z)-B_1)$ is contained in $\nsupp(v-B\cdot(y+z))$.

This containment might not be strict, but certainly
$(v-B\cdot(y+z)-B_1)_3 > (v-B\cdot(y+z))_3$ (or
$(v-B\cdot(y+z)-B_2)_3 > (v-B\cdot(y+z))_3$ in the other cases).
Moreover, we can repeat this process, and keep adding columns of $B$ until
the third coordinate is a nonnegative integer, while keeping
the first, second and fourth coordinates also nonnegative.
This shows that $u-B\cdot z = v-B\cdot(y+z)$ does not have
minimum negative support, which is the contradiction we wanted.
\end{proof}

Now we can look at the
logarithm-free canonical series solution $\phi_{u-e_3}$
of $H_A(\beta)$ corresponding to the fake exponent $u-e_3$. 
We claim that this function does not lie in the image of the map $\partial_3$
between the solution spaces of $H_A(A\cdot v)$ and $H_A(\beta)$.

\begin{proposition}
\label{hardone}
If $\psi$ is a solution of $H_A(A\cdot v)$ and $u$ is as in
Lemma \ref{littlebug},
then $\partial_3 \psi \neq \phi_{u-e_3}$.
\end{proposition}

\begin{proof}
Suppose there is a solution $\psi$ of $H_A(A\cdot v)$ such that 
$\partial_3 \psi = \phi_{u-e_3}$. We will obtain a contradiction.

We proceed as in the part of the proof of 
Theorem \ref{kernel} where we show that the functions $\varphi_i$
are logarithm-free. The first step is to use
the Observation \ref{keyobservation} to write
$\psi=\psi_{\delta}+c_{\delta}x^u \log(x)^{\delta}$
for every $\delta \in {\mathcal{S}}_{\max}$.
We apply Lemma \ref{bigbug}, with the goal
of showing that $\psi$ has no terms in $\log(x_j)$ for
$j \in {\mathcal{I}}$. Let $j \in {\mathcal{I}}$,
pick $\delta \in {\mathcal{S}}_{\max}$ maximal with
respect to the $j$-th coordinate, and let $z=z^{(j)}$
from Lemma \ref{bigbug}. Then
\[ \partial^{(B\cdot z)_-} \psi = \partial^{(B\cdot z)_+} 
\psi_{\delta} + 
\partial^{(B\cdot z)_+} c_{\delta} x^u \log(x)^{\delta} \; . \]
As in the proof of Theorem \ref{kernel}, there are nonzero
terms when we compute
$\partial^{(B\cdot z)_+} c_{\delta} x^u \log(x)^{\delta}$
using the product rule.
Now, all these terms 
have either logarithms or denominator a positive integer 
power of $x_j$. 

By construction, $u_j = \eta_j \geq 0$, so that $j \not \in \nsupp(u-e_3)$.
Now, since $(B\cdot z)_3 < 0$, $\partial^{(B\cdot z)_-} \psi$
is a further derivative of $\phi_{u-e_3}$.
But
$\phi_{u-e_3}$ has no terms with denominator $x_j^k$ with $0< k  \in \N$.
This means that 
$\partial^{(B\cdot z)_+} c_{\delta} x^u \log(x)^{\delta}$
must be cancelled with terms coming from
$\partial^{(B\cdot z)_+} \psi_{\delta}$, and this implies
(again, as in Theorem \ref{kernel}) that $\delta_j=0$ or,
equivalently, that $\psi$ has no terms in $\log(x_j)$
for $j \in {\mathcal{I}}$. From this we can
show that $\psi$ is logarithm-free.

Now, $\phi_{u-e_3}$ has a term $x^{u-e_3}$, and the
only way this term matches with a term of $\partial_3 \psi$ is
if $\psi$ has a term $x^u \log(x_3)$. But $\psi$ is logarithm-free,
and we obtain a contradiction.
\end{proof}

It is now time to deal with logarithmic solutions
of $H_A(A\cdot v)$ corresponding to 
exponents that differ by an integer vector with $v$.

\begin{proposition}
\label{prelaurentimage}
If $\psi$ is a solution of $H_A(A\cdot v)$ such that $\partial_3 \psi$
lies in $\mbox{Span}\, \{ \phi_1,\phi_2,\phi_3,\phi_4\}$,
where the functions $\phi_i$ are the solutions of $H_A(\beta)$
we introduced in Lemma \ref{pseudolaurentsolutions}, 
then (modulo the kernel of $\partial_3$),
\[ \psi = \tilde{\psi} + x^v \sum_{i=1}^n c_i \log(x_i) \; ,\]
where $\tilde{\psi}$ is a logarithm free series with exponents that differ
by integer vectors with $v$, 
that has no
term in $x^v$, and the vector $(c_1,\dots,c_n)^t$ belongs to the kernel of
$A$. 
\end{proposition}

\begin{proof}
Pick any solution $\psi$ of $H_A(A\cdot v)$ whose 
derivative  with respect
to $x_3$ lies in $\mbox{Span}\, \{ \phi_1,\phi_2,\phi_3,\phi_4\}$.
Write $\psi$ as in Observation \ref{keyobservation}:
\[ \psi =\psi_{\delta} +  \log(x)^{\delta} f \, , \]
for $\delta \in {\mathcal{S}}_{\max}$.

Here we must have $f = c_{\delta} x^v$, since $\partial_3 \psi$
lies in $\mbox{Span}\, \{ \phi_1,\phi_2,\phi_3,\phi_4\}$.
Suppose that $\delta_3 \neq 0$.
Look at the logarithm-free function
\[  \partial_3\, \psi = 
\partial_3\, \psi_{\delta} +c_{\delta} \log(x)^{\delta-e_3} 
\frac{x^v}{x_3} \; .\]
If $\delta \neq e_3$,
the term $c_{\delta} \log(x)^{\delta-e_3} \frac{x^v}{x_3}$ has
logarithms, so it must be cancelled with terms coming from
$\partial_3 \, \psi_{\delta}$. Thus 
$\psi_{\delta}$ must have a sub-series $g$ such that $\partial_3 g = 
c_{\delta} \log(x)^{\delta-e_3} \frac{x^v}{x_3}$. Then
$g-c_{\delta} x^v \log(x)^{\delta}$ is constant with respect to $x_3$,
which contradicts the construction of $\psi_{\delta}$ (all its logarithmic 
terms are either less than or incomparable to $\delta$).
Therefore, $\delta_3=0$ or $\delta=e_3$.

Now choose $\delta \in {\mathcal{S}}_{\max}$
with $\delta_2$ maximal. Suppose $\delta_2 \geq 1$. Consider the
function:
\begin{eqnarray*}
\partial^{(B_2)_-} \psi &=& \partial^{(B_2)_+} \psi \\
   & = &  \partial^{(B_2)_+} \psi_{\delta} + c_{\delta} \log(x)^{\delta-e_2} 
\frac{\partial^{(B_2)_+-e_2}x^v}{x_2}+ \\
   & & +\; \mbox{other terms coming from $c_{\delta}\partial^{(B_2)_+} x^v 
\log(x)^{\delta}$}
\end{eqnarray*}
If $\delta \neq e_2$ the nonzero summand 
$c_{\delta} \log(x)^{\delta-e_2}\frac{\partial^{(B_2)_+-e_2}x^v}{x_2}$
has logarithms, hence it must be cancelled by some other term of 
the right hand side sum. Since the numerator 
$\partial^{(B_2)_+-e_2}x^v$ is constant with respect to 
$x_2$, this is impossible.
Thus $\delta=e_2$ or $\delta_2=0$.

Similar arguments using $B_1+B_2$ show that $\delta \in {\mathcal{S}}_{\max}$
maximal with respect to the first coordinate must be either $e_1$
or have $\delta_1=0$, and the same using $B_1$ will give the
analogous conclusion when $\delta$ is maximal with respect to
the fourth coordinate. 

Now pick $i>4$ and choose $\delta \in {\mathcal{S}}_{\max}$ with $\delta_i$ 
maximal.
Suppose
$\delta_i>1$. Then $\delta_l=0$ for $1 \leq l \leq 4$.

We know (looking at the homogeneities (\ref{homogeneities})) that
\[ (A\cdot v)_j \psi 
= \sum_{k=1}^n a_{jk} \theta_k \psi \;, \;\; j=1,\ldots,d \, .\]
Call $\tilde{c}$ the coefficient of $\log(x)^{\delta-e_i}x^v$ in $\psi$.
Comparing both sides of the previous equalities, we conclude that
\[ \tilde{c} (A\cdot v)_j = \sum_{\gamma \in {\mathcal{S}}_{\max} : 
\gamma-e_k = \delta-e_i}
a_{jk} c_{\gamma} \;, \;\; j=1,\ldots,d \, .\]
Let $v'$ be the vector whose $k$-th coordinate is $c_\gamma$ if $\gamma-e_k = 
\delta-e_i$, and the rest are zeros. If $\tilde{c}=0$, $v'$ is a nonzero 
element
of the kernel of $A$, whose first $4$ coordinates are $0$. Such an element
does not exist. If $\tilde{c} \neq 0$, $A(1/\tilde{c})v'=A\cdot v$, 
and the first
$4$ coordinates of $(1/\tilde{c})v'$ are zero. Thus, $v-(1/\tilde{c})v'$
is a vector in the kernel of $A$ and
is therefore of the form $B\cdot z$, 
for $z\in \R^2$. Since the first $4$ rows of
$B$ lie in different quadrants of $\Z^2$
and the first four entries of $v-(1/\tilde{c})v'$
are nonnegative, this is impossible.

Hence
\[ \psi = \tilde{\psi} + x^v \sum_{i=1}^n c_i \log(x_i) \; , \]
where $\tilde{\psi}$ is logarithm-free. If we assume that $\tilde{\psi}$
has no term $x^v$ (perfectly legal, since this is a solution
of $H_A(A\cdot v)$ that is constant with respect to $x_3$), the fact
that the vector $(c_1,\dots ,c_n)^t$ belongs to the 
kernel of $A$ follows from the homogeneities (\ref{homogeneities}).
\end{proof}

\begin{theorem}
\label{Laurentimage}
The dimension of the intersection of the image of the map $\partial_3$
and the space
$\mbox{Span}\, \{ \phi_1,\phi_2,\phi_3,\phi_4\}$ 
is at most $2$. Consequently,
the dimension of $\mbox{Span}\, \{ \phi_1,\phi_2,\phi_3,\phi_4\}$
modulo the image of $\partial_3$ is greater than or equal to $2$,
that is, we can choose two functions in 
$\mbox{Span}\, \{ \phi_1,\phi_2,\phi_3,\phi_4\}$ which span
a $2$-dimensional subspace of the cokernel of $\partial_3$.
\end{theorem}

\begin{proof}
By Proposition \ref{prelaurentimage}, an element $\psi$
of the solution space of $H_A(A\cdot v)$ 
such that $\partial_3 \psi$ 
lies in $\mbox{Span}\, \{ \phi_1,\phi_2,\phi_3,\phi_4\}$
is of the form
\[ \psi = \tilde{\psi} + x^{v} \sum_{i=1}^n c_i \log(x_i)\; , \]
with $\tilde{\psi}$ a logarithm-free function with integer 
exponents, no term $x^v$,
and $(c_1,\dots , c_n)^t$ in the kernel of $A$.

Notice that once the $c_i$ are fixed, $\psi$ is unique with those $c_i$,
since the difference of two such functions would be a logarithm-free
solution of $H_A(A\cdot v)$ with no term in the kernel of $\partial_3$,
whose derivative with respect to $x_3$ belongs to 
$\mbox{Span}\, \{ \phi_1,\phi_2,\phi_3,\phi_4\}$.
It follows
from Proposition \ref{monokernel} that this difference must be zero.

Since the vector of the $c_i$ is in the kernel of $A$, the previous
remark implies that
the space of solutions of $H_A(A\cdot v)$ 
whose derivative with respect 
to $x_3$ lies in
$\mbox{Span}\, \{ \phi_1,\phi_2,\phi_3,\phi_4\}$
has dimension at most $2$, the
dimension of the kernel of $A$.
\end{proof}

\begin{lemma}
\label{lind}
The functions $\phi_{u-e_3}$ constructed in Proposition \ref{hardone} 
and the functions 
from Theorem \ref{Laurentimage} that span a $2$-dimensional
subspace of the cokernel of $\partial_3$
are linearly independent modulo the
image of the map $\partial_1$.
\end{lemma}

\begin{proof}
By contradiction, suppose there is a solution $\psi$ of $H_A(A\cdot v)$ 
such that
\[ \partial_1 \psi = L + \left(
\sum_{u \not \in \Z^n : x^u \in \ker(\partial_3)}  c_u \phi_{u-e_3}  
\right) \; ,\]
where $L$ is a linear combination of the functions from Theorem
\ref{Laurentimage}. We can write 
\[ \psi= \psi_L + \left(
\sum_{u \not \in \Z^n : x^u \in \ker(\partial_3)} \psi_u  \right) \; ,\]
where $\psi_u$ is the sum of the terms in $\psi$ whose exponents and $u$
differ by an integer vector, and $\psi_L$ is the sum of the 
terms in $\psi$ whose
exponents and $v$ differ by an integer vector. 
Clearly, $\partial_3 \psi_u = c_u \phi_{u-e_3} $
and $\partial_3 \psi_L=L$. But the functions $\psi_u$ and $\psi_L$
must be solutions of $H_A(A\cdot v)$. To see this, notice that
if two monomials $x^{\alpha_1}$ and $x^{\alpha_2}$ are such that
$\alpha_1-\alpha_2 \not \in \Z^n$, then the intersection of the $D$-modules
obtained by acting with $D$ on 
$x^{\alpha_1}$ and $x^{\alpha_2}$ is either
empty or $\{ 0 \}$.

Therefore all the $c_u$ must be zero (by Proposition
\ref{hardone}) and also $L$ must be zero (by Theorem \ref{Laurentimage}).
\end{proof}

We have now all the ingredients to show that the 
parameter $\beta$ from Construction 
\ref{construction} is indeed an exceptional parameter.

\begin{theorem}
\label{rankjump}
Let $\beta$ be the parameter from Construction \ref{construction}. Then 
\[ \rank(H_A(\beta)) \geq \vol(A) + 1 \, .\]
\end{theorem}

\begin{proof}
In Proposition \ref{hardone} and Theorem \ref{Laurentimage}
we built one function in $\coker(\partial_3)$ for 
each function in a basis of $\ker(\partial_3)$ (which we
knew from Theorem \ref{kernel}). Moreover,
Theorem \ref{Laurentimage} provided at least two linearly
independent functions for $x^v$. Lemma \ref{lind} shows that
all of these functions are linearly independent. Therefore
\[ \dim(\coker(\partial_1)) \geq \dim(\ker(\partial_3)) +1 \, ,\]
and this implies that:
\begin{equation}
\label{finally}
\dim(\coker(\partial_1))+ \dim(\mbox{im}(\partial_3)) 
\geq \dim(\ker(\partial_3))+ \dim(\mbox{im}(\partial_3)) +1 \, ,
\end{equation}
where $\mbox{im}(\partial_3)$ is the image of $\partial_3$.
The left hand side of (\ref{finally}) equals the dimension
of the solution space of $H_A(\beta)$. The right hand side equals
$1$ plus the dimension of the solution space of $H_A(A\cdot v)$.
This concludes the proof.
\end{proof}

When $n>4$, we can use Theorem \ref{rankjump} to reach a stronger
conclusion about the exceptional set of $A$.

\begin{theorem}
\label{exceptionalset}
Let $A$ be such that $I_A$ is a non Cohen Macaulay toric ideal,
with a Gale diagram whose first four rows meet each of the 
open quadrants of $\Z^2$. Let $v_1$, $v_2$, $v_4$ be as in
Construction \ref{construction}.
If $n > 4$, the $(n-4)$-dimensional affine space
parametrized by:
\[ (s_5,\dots,s_n) \longmapsto
A \cdot (v_1 e_1+v_2 e_2 -e_3 + v_4 e_4 + \sum_{i=5}^n s_i e_i)\]
is contained in the exceptional set ${\mathcal{E}}(A)$
In particular, ${\mathcal{E}}(A)$ is an infinite set.
\end{theorem}

\begin{proof}
Pick $(s_5,\dots,s_n) \in \C^{n-4}$, and 
$\alpha_5,\dots,\alpha_n$ as in Construction \ref{construction}.
We can choose $\kappa_0$ small enough so that the
numbers $\tilde{\alpha}_i = s_i + \kappa \alpha_i $
satisfy the conditions of Construction \ref{construction} for
all $0 < \kappa < \kappa_0$. Call
\[ \beta_{\kappa}:= A \cdot (v_1 e_1+v_2 e_2 -e_3 + v_4 e_4 
+ \sum_{i=5}^n \tilde{\alpha}_i e_i) \, ,\]
and 
\[ \beta:= A \cdot (v_1 e_1+v_2 e_2 -e_3 + v_4 e_4 
+ \sum_{i=5}^n s_i e_i) \, .\]
Then Theorem \ref{rankjump} implies that 
$\rank(H_A(\beta_{\kappa})) \geq \vol(A) + 1$
for all $0 < \kappa < \kappa_0$. Now the proof of Theorem 3.5.1
in \cite{SST} implies that $\rank(H_A(\beta)) \geq \vol(A)+1$.
This concludes the proof.
\end{proof}

\section{Examples and Final Remarks}
\label{fifthsection}
To conclude, we illustrate our ideas in an example, and point out 
some open questions about rank jumps, even in codimension 2. We choose
the following matrix:
\[ A = \left( \begin{array}{ccccc} 1 & 1 & 1 & 1 & 1 \\
0 & 1 & 0 & 1 & 0 \\ 0& 0& 1 & 1 & -2 \end{array} \right)\]
In this case, $\vol(A)=4$. Since $A$ has the Gale diagram
\[B=\left( \begin{array}{cc} 1 & 2 \\ -1 & 1 \\ -1 & -1 \\ 1 &-1 \\ 
0 & -1 \end{array} \right) \]
which meets the four open quadrants of $\Z^2$, we conclude that the toric
ideal $I_A$ is not Cohen Macaulay. 

We will show that 
\[ {\mathcal{E}}(A)=\{ (1,0,-1)^t + s (1,0,-2)^t:
s \in \C \} \, . \]
From Theorem \ref{exceptionalset} we conclude that
the line $\{ A \cdot (2,0,-1,0,s)^t : s \in \C\}=\{ (1,0,-1)^t + s (1,0,-2)^t:
s \in \C \}$
is contained in the exceptional set ${\mathcal{E}}(A)$.

In Section 4.6 of \cite{SST}, we 
see that, for each initial monomial ideal of $I_A$, we can construct
a finite arrangement of planes in $\C^d$ that contains the exceptional set.
It is therefore informative to compute all initial monomial ideals
of $I_A$, form the corresponding arrangement for each initial ideal,
and intersect all of them. In our example, $I_A$ has $9$ initial monomial
ideals (computed using {\tt TiGERS}, \cite{tigers}). 
The intersection of all the arrangements coming from these initial ideals
is the zero set of the ideal:
\[ \langle \beta_2, 2\beta_1+\beta_3-1 \rangle \cap \langle \beta_3-1,
\beta_2-3,\beta_1-4\rangle \cap \]
\[\cap \langle \beta_3-1, \beta_2-1, \beta_1-1 \rangle
\cap \langle \beta_3^2, \beta_2-3\beta_3, \beta_1 -4\beta_3 \rangle \]
that is, the union of the points:
\[ (0,0,0)^t, (4,3,1)^t,(1,1,1)^t \]
and the line:
\[ \{ (1,0,-1)^t + s (1,0,-2)^t : s \in \C \} \, . \]
With the help of {\tt Macaulay2} for the Weyl Algebra (which can
be obtained at \cite{m2}), we
find that the points $(0,0,0)^t$ , $ (4,3,1)^t$, $(1,1,1)^t$ 
do not belong to ${\mathcal{E}}(A)$, so 
that $ \{ (1,0,-1)^t + s (1,0,-2)^t : s \in \C \} \supseteq 
{\mathcal{E}}(A) \, $.
We conclude that
\[\{ (1,0,-1)^t + s (1,0,-2)^t : s \in \C \} = 
{\mathcal{E}}(A) \, .\]

It is an open question to
give a sharp bound for the maximum possible magnitude 
of rank jumps, even in codimension 2. Corollary 4.1.2 in \cite{SST}
gives the only known upper bound for the rank of 
$A$-hypergeometric systems, but most
likely, it is far from optimal. 
Also, we can find examples in codimension $2$ of exceptional parameters 
where the rank jump is more than $1$. For instance, 
let 
\[ A=\left( \begin{array}{cccccc} 1 & 1 & 1 & 1 & 1 & 1  \\
0 & 0 & 0 & 0 & -2 & 1 \\
1 & 0 & 0 & 1 & 0 & 0 \\
1 & 0 & 2 & 0 & 1 & 1 \end{array} \right) \; .\]
Then we have $\vol(A)=9$. Computing a Gale diagram, we see that $I_A$ is not
Cohen Macaulay.
In fact, Theorem \ref{exceptionalset} produces, for instance, $\beta=(4,2,0,5)^t$, with
$\rank(H_A(\beta))=10$. However,
for $\beta=(2,1,0,2)^t$, $\rank(H_A(\beta))=11$. This parameter vector
also comes from Theorem \ref{exceptionalset}.

\medskip

There is hope that the construction in this article can be 
extended to provide exceptional parameters for $A$-hypergeometric
systems such that certain reverse lexicographic initial ideals of $I_A$
have embedded primes. Work in that direction is ongoing.

\medskip

\noindent {\bf Acknowledgments:}
I am very grateful to Bernd Sturmfels: this work grew out of our weekly
discussions. I would also like to thank Rekha Thomas, who told me
about the wonderful polytopes $P^{\bar{\sigma}}_{\eta}(0)$, and
Harrison Tsai, for listening about this project in its various stages of
completion, and for all the {\tt Macaulay2} help.

\end{document}